\documentclass[letterpaper, fleqn, 11pt]{amsart}
\usepackage[english]{babel}
\usepackage{graphicx,color}
\usepackage[all]{xy}
\usepackage{amsfonts,dsfont,mathrsfs}
\usepackage{amsmath, amsthm, amssymb}
\usepackage{enumerate, url}

\newcommand{\bbR}{\ensuremath{\mathbb{R}}}

\newcommand{\bbZ}{\ensuremath{\mathbb{Z}}}

\newcommand{\bbH}{\ensuremath{\mathbb{H}}}

\DeclareMathOperator{\id}{id}

\DeclareMathOperator{\rk}{rk}
\DeclareMathOperator{\diam}{diam}
\DeclareMathOperator{\injrad}{injrad}
\DeclareMathOperator{\vol}{vol}
\DeclareMathOperator{\dvol}{dvol}
\DeclareMathOperator{\Ric}{Ric}
\DeclareMathOperator{\isom}{Isom}
\DeclareMathOperator{\comm}{Comm}

\DeclareMathOperator{\Out}{Out}

\numberwithin{equation}{section}

\newtheorem{thmnr}{Theorem}[section]

\newtheorem{propnr}[thmnr]{Proposition}

\newtheorem{lemnr}[thmnr]{Lemma}

\newtheorem{cornr}[thmnr]{Corollary}

\newtheorem{conjnr}[thmnr]{Conjecture}
\theoremstyle{definition}

\newtheorem{dfnnr}[thmnr]{Definition}

\newtheorem{rmknr}[thmnr]{Remark}

\newtheorem{claimnr}[thmnr]{Claim}

\begin{document}

\title{Symmetry gaps in Riemannian geometry and minimal orbifolds}

\author{Wouter van Limbeek}

\date{\today}

\begin{abstract} We study the size of the isometry group Isom$(M,g)$ of Riemannian manifolds $(M,g)$ as $g$ varies. For $M$ not admitting a circle action, we show that the order of Isom$(M,g)$ can be universally bounded in terms of the bounds on Ricci curvature, diameter, and injectivity radius of $M$. This generalizes results known for negative Ricci curvature to all manifolds.

 More generally we establish a similar universal bound on the index of the deck group $\pi_1(M)$ in the isometry group $\isom(\widetilde{M},\widetilde{g})$ of the universal cover $\widetilde{M}$ in the absence of suitable actions by connected groups. We apply this to characterize locally symmetric spaces by their symmetry in covers. This proves a conjecture of Farb and Weinberger with the additional assumption of bounds on curvature, diameter, and injectivity radius. Further we generalize results of Kazhdan-Margulis and Gromov on minimal orbifolds of nonpositively curved manifolds to arbitrary manifolds with only a purely topological assumption.
\end{abstract}

\maketitle

\tableofcontents

\section{Introduction}
\label{sec:intr}
Let $M$ be a closed Riemannian manifold with metric $g$. A general question is to study the relation between the geometry of $M$ and the group of isometries $\isom(M,g)$. An important example is the classical result of Bochner-Yano \cite{bochneryano} that if $M$ is negatively curved then $\isom(M)$ is finite.

Given this qualitative result, one can hope for a quantitative version that bounds the size of the isometry group. The first result in this direction is a classical theorem due to Hurwitz, which states that for any hyperbolic metric on a surface $\Sigma$ of genus $g\geq 2$, the order of $\isom(\Sigma)$ is at most $84(g-1)$. This result already shows there is an interesting connection between the question of bounds on the order of Isom$(M,g)$ and the topology of $M$. Further results in this direction have been proved by Huber \cite{isomhyp} (for hyperbolic manifolds), Im Hof (for negatively curved manifolds) \cite{isomneg}, Maeda (for nonpositively curved manifolds with negative Ricci curvature) \cite{isomnonpos}, and Katsuda (for manifolds with negative Ricci curvature) \cite{isomnegric}, and Dai-Shen-Wei \cite{DSW}. For more information see \cite{DSW} and the references therein. We will just state the latter result since it is the most general.
	\begin{thmnr}[{\cite{DSW}}] Let $n\geq 1$ and $\Lambda>\lambda, \varepsilon_0, D$ be positive constants. Then there exists $C\geq 1$ such that for any closed Riemannian $n$-manifold $(M,g)$ such that
	$$-\Lambda\leq \Ric_M \leq -\lambda<0, \hspace{0.5 cm} \injrad_M \geq \varepsilon_0, \hspace{0.5 cm} \diam(M)\leq D,$$
we have $|\isom(M)|\leq C$. 
	\label{thm:dsw}
	\end{thmnr}
The proofs of the above results rely on differential geometric techniques in negative (Ricci) curvature such as Bochner identities. On the other hand, it is clear that the Bochner-Yano theorem generalizes to many manifolds that are not negatively curved. For instance, Borel proved that if $M$ is closed aspherical and $\pi_1(M)$ is centerless then $M$ admits no nontrivial action by a compact connected Lie group \cite{borelnoact}. Hence for any metric $g$ on $M$, the isometry group $\isom(M,g)$ is finite. 

\subsection*{Main results.} As far as we are aware no progress has been made on quantitative generalizations of the Bochner-Yano theorem without some assumption of negative curvature. The goal of this paper is to prove results similar to Theorem \ref{thm:dsw} outside of the setting of negative curvature. However, $M$ may admit metrics with infinite isometry group so that no upper bound on the number of isometries can be obtained. Our main result is that the existence of an action by a connected Lie group is the only obstruction to a uniform bound on the size of isometry groups. More precisely, we prove the following.
	\begin{thmnr} Let $\Lambda, \varepsilon, D>0$ and $n\geq 1$. Then there exists $C>0$ with the following property. Suppose $(M,g)$ is a closed Riemannian $n$-manifold such that $M$ does not admit a $C^2$ circle action, and such that
			\begin{equation}
			|\Ric_g|\leq \Lambda, \hspace{0.5 cm} \injrad(M,g)\geq \varepsilon, \hspace{0.5 cm} \diam(M,g)\leq D.
			\label{eq:bounds}
			\end{equation}
	Then $\isom(M,g)$ has order at most $C$.
	\label{thm:cmptbddbound}
	\end{thmnr}
	\begin{rmknr} At least some of the bounds on Ricci curvature, injectivity radius and diameter are essential in Theorem \ref{thm:cmptbddbound} because of the following example. Cappell-Weinberger-Yan \cite{connerraymondconj} constructed smooth closed aspherical manifolds $M$ of any dimension $>7$ that admit no continuous action by a connected group. However, as was pointed out to me by Shmuel Weinberger, these examples have a degree 3 self-cover $p:M\rightarrow M$. From the construction it is clear that $p$ and all of its iterates $p^n:M\rightarrow M$ are regular covers, so that $M$ admits a smooth effective action by a group of order $3^n$ for any $n\geq 1$.\label{rmk:excsym}\end{rmknr}
	\begin{rmknr} As explained above, if $M$ is a manifold with a $C^\infty$ circle action, then there is a smooth metric such that this action is by isometries. In Theorem \ref{thm:cmptbddbound} we assume that there is no $C^2$ action. We do not know if the nonexistence of a $C^\infty$ circle action suffices to obtain the conclusion of Theorem \ref{thm:cmptbddbound}.\end{rmknr}
	\begin{rmknr} Our proof gives no information about the size of $C$, but it would be interesting to obtain explicit upper bounds.\end{rmknr}
As mentioned above, in the special case that $M$ is a closed aspherical manifold and $\pi_1(M)$ is centerless, Borel proved that no connected Lie group acts effectively on $M$. Therefore Theorem \ref{thm:cmptbddbound} yields a quantitative version of Borel's theorem. In this case we actually obtain an explicit, but probably very far from optimal, bound on $C$.
	\begin{thmnr}[Quantitative Borel] Let $\Lambda, \varepsilon, D>0$ and $n\geq 1$. Then there exists $C>0$ with the following property. Let $(M,g)$ be a closed aspherical Riemannian $n$-manifold with centerless fundamental group and satisfying the bounds of Equation \ref{eq:bounds}. Then $\isom(M,g)$ has order at most $C$.
	
Set $N=\frac{V_k\left(\frac{D}{2}\right)}{V_k(\frac{\varepsilon}{4})}$, where $V_k(r)$ denotes the volume of a ball of radius $r$ in a simply-connected manifold of dimension $n$ and constant curvature $\frac{-\Lambda}{n-1}$. Then we can take $C=N^N$.
	\label{thm:quantborel}
	\end{thmnr}
Note that in Theorems \ref{thm:cmptbddbound} and \ref{thm:quantborel} we do not require the Ricci curvature to be negative. Consequently, the proofs will not use phenomena in negative curvature and hence will be very different from the proof of Theorem \ref{thm:dsw} (and the previous results of Bochner-Yano, H\"uber, Im Hof, Katsuda, Maeda, and Dai-Shen-Wei). Instead our methods are more topological. We will give a further outline of the proof of Theorem \ref{thm:cmptbddbound} at the end of this section.

The use of topological tools (instead of geometry in negative curvature) has the advantage of generalizing well to covers, and this yields far more information. The viewpoint of studying the isometries of covers of $M$, rather than just of $M$ itself, has been first considered by Eberlein \cite{eblatt, ebisom} (for nonpositively curved manifolds) and later Frankel \cite{frharm} (for semisimple isometry groups), and since greatly developed by Farb-Weinberger \cite{FWarith,FW} (for aspherical manifolds). Here one studies the relation between the geometry of $(M,g)$ and the isometry group $\isom(\widetilde{M})$ of the universal cover $\widetilde{M}$. Note that the group $I(\widetilde{M})$ potentially contains much more information than $\isom(M)$. For example, if $M$ is hyperbolic, then the Bochner-Yano theorem implies that $\isom(M)$ is finite, but $\isom(\widetilde{M})$ is a Lie group that acts transitively on $\widetilde{M}$.
 
Of course $\isom(\widetilde{M})$ always contains the deck group $\pi_1(M)$. In \cite{FW}, Farb-Weinberger prove that if $M$ is aspherical and $\isom(\widetilde{M})$ contains the deck group $\pi_1(M)$ with infinite index, then $M$ is a Riemannian orbibundle with locally homogeneous fibers (see Section \ref{sec:fw} for more information and explanation of this terminology). 

The correct analogue of the order of $\isom(M,g)$ is the index $[\isom(\widetilde{M},\widetilde{g}):\pi_1(M)]$. Note that the Bochner-Yano theorem does not hold for covers, because, as remarked above, if $M$ is hyperbolic then $[\isom(\widetilde{M}):\pi_1(M)]=\infty$. Still, the analogue of Theorem \ref{thm:cmptbddbound} is true. Namely we bound the index $[\isom(\widetilde{M},\widetilde{g}):\pi_1(M)]$ in terms of $\dim M$, bounds on Ricci curvature, diameter, and injectivity radius of $M$ under the assumption of absence of appropriate actions.
	\begin{thmnr} Let $\Lambda, \varepsilon, D>0$ and $n\geq 1$. Then there exists $C>0$ with the following property. Suppose $(M,g)$ is a closed Riemannian $n$-manifold such that $\widetilde{M}$ does not admit a proper $C^2$ action by a postive-dimensional Lie group $G$ containing $\pi_1(M)$, and such that
			\begin{equation}
			|\Ric_g|\leq \Lambda, \hspace{0.5 cm} \injrad(M,g)\geq \varepsilon, \hspace{0.5 cm} \diam(M,g)\leq D.
			\end{equation}
	Then $[\isom(\widetilde{M}):\pi_1(M)]\leq C$. 
	\label{thm:coverbddbound}
	\end{thmnr}
We prove a stronger result in the special case of aspherical manifolds $M$ such that $\pi_1(M)$ contains no nontrivial normal abelian subgroup. Namely, in this case we are able to prove that the above result holds with $C^2$ actions replaced by $C^\infty$ actions. Therefore we obtain the existence of a gap in smooth symmetry of these manifolds. Namely, either $[\isom(\widetilde{M}):\pi_1(M)]\leq C$ or $[\isom(\widetilde{M}):\pi_1(M)]=\infty$. Combined with the work of Farb-Weinberger on the latter case \cite{FW}, this yields the following theorem.
	\begin{thmnr} Let $\Lambda, \varepsilon, D>0$ and $n\geq 1$. Then there exist $C, d>0$ with the following property. Let $(M,g)$ be a closed aspherical Riemannian $n$-manifold such that $\pi_1(M)$ contains no nontrivial normal abelian subgroup, and such that
		\begin{equation}
			|\Ric_g|\leq \Lambda, \hspace{0.5 cm} \injrad(M,g)\geq \varepsilon, \hspace{0.5 cm} \diam(M,g)\leq D.
		\end{equation}
	Then $[\isom(\widetilde{M}):\pi_1(M)]\leq C$ or $M$ has a finite cover $M'$ of degree at most $d$ such that $M'$ is isometric to a nontrivial Riemannian warped product $B\times_f N$ where $f:B\rightarrow\bbR_{>0}$ and for every $b\in B$, the copy $\{b\}\times N$ of $N$ is isometric to a locally symmetric space of noncompact type.
	\label{thm:coverreg}
	\end{thmnr}
	
\subsection*{Applications.} We now give two applications of Theorem \ref{thm:coverreg}. The first application characterizes locally symmetric spaces in terms of isometries of the universal cover. In this context, Farb-Weinberger proved the following.
	\begin{thmnr}[Farb-Weinberger {\cite[Theorem 1.3]{FW}}] Let $M$ be a closed aspherical, smoothly irreducible, Riemannian manifold such that $\pi_1(M)$ contains no nontrivial normal abelian subgroups and $[\isom(\widetilde{M}):\pi_1(M)]=\infty$. Then $M$ is isometric to a locally symmetric space of noncompact type.
	\label{thm:charlocsym}
	\end{thmnr}
Here a manifold $M$ is called \emph{smoothly irreducible} if there is no finite cover of $M$ that is diffeomorphic to a nontrivial product. In addition Farb-Weinberger conjectured the following quantitative version of Theorem \ref{thm:charlocsym}.
	\begin{conjnr}[Farb-Weinberger {\cite[Conjecture 1.6]{FW}}] Let $M$ be a smooth closed, aspherical, smoothly irreducible manifold such that $\pi_1(M)$ contains no nontrivial normal abelian subgroups. Then there exists $C>1$ only depending on $\pi_1(M)$ such that if $g$ is any Riemannian metric on $M$ with $[\isom(\widetilde{M}):\pi_1(M)]\geq C$, then $M$ is isometric to a locally symmetric space of noncompact type.
	\label{conj:magicnr}
	\end{conjnr}
	\begin{rmknr} The assumption that $\pi_1(M)$ has no nontrivial normal abelian subgroups can be explained as follows. First note that the conjecture fails for tori. More generally there are many fiber bundles $M\rightarrow N$ with fibers isometric to tori for which the conjecture fails. However, in this case $M$ has a nontrivial normal abelian subgroup.
	
	On the other hand suppose $M$ is closed aspherical and $\pi_1(M)$ has trivial center. Then by a theorem of Borel, $M$ admits no effective actions by connected Lie groups.\end{rmknr}
	\begin{rmknr} In view of Remark \ref{rmk:excsym}, the assumption in Conjecture \ref{conj:magicnr} that $\pi_1(M)$ has no nontrivial normal abelian subgroups cannot be weakened to assuming that $M$ admits no action by a connected Lie group.\end{rmknr}
As evidence for Conjecture \ref{conj:magicnr}, Farb-Weinberger proved the conjecture if $M$ is assumed to be diffeomorphic to a locally symmetric space of noncompact type. Avramidi \cite{grigoriperflat} proved the conjecture for noncompact finite volume locally symmetric spaces.
Tam Nguyen Phan \cite{tamminorb} proved Conjecture \ref{conj:magicnr} for piecewise locally symmetric manifolds. These manifolds are obtained by gluing noncompact finite volume locally symmetric spaces.

We will prove Conjecture \ref{conj:magicnr} for general $M$, but subject to bounds on the Ricci curvature, injectivity radius, and diameter of $(M,g)$. The following is immediate from Theorem \ref{thm:coverreg} and the assumption of smooth irreducibility of $M$.
	\begin{cornr} Let $\Lambda, \varepsilon, D>0$ and $n\geq 2$. Then there exists $C>0$ with the following property. Suppose $(M,g)$ is a closed, aspherical, smoothly irreducible Riemannian $n$-manifold such that $\pi_1(M)$ contains no nontrivial normal abelian subgroup,  and
		\begin{equation}
		|\Ric_g|\leq \Lambda, \hspace{0.5 cm} \injrad(M,g)\geq \varepsilon, \hspace{0.5 cm} \diam(M,g)\leq D.
		\end{equation}
	Then $[\isom(\widetilde{M}):\pi_1(M)]\leq C$ or $(M,g)$ is isometric to a locally symmetric space of noncompact type.
	\label{cor:bddmagicnr}
	\end{cornr}
The second application of Theorem \ref{thm:coverbddbound} is to the phenomenon of minimal orbifolds. This was first discovered in the context of symmetric spaces by Kazhdan-Margulis \cite{minorbhomog}. They proved that for every $n\geq 1$, the volume of locally symmetric orbifolds of dimension $n$ is bounded below by some number $\mu_n>0$. More precisely, if $X$ is a symmetric space of noncompact type of dimension $n$, then for any group $\Gamma$ acting on $X$ properly discontinuously (not necessarily freely), we have $\vol(X\slash\Gamma)\geq \mu_n$. This was generalized to manifolds with negative curvature by Gromov \cite{gromovneg}. We prove the following related result on contractible complete Riemannian manifolds $X$ that does not assume negative curvature.
	\begin{cornr}[Minimal orbifolds] Let $\Lambda, \varepsilon, D>0$ and $n\geq 2$. Then there exists $\mu>0$ with the following property. Let $X$ be a contractible Riemannian $n$-manifold with $|\Ric|\leq\Lambda$ and admitting a compact manifold quotient $M$ such that $\pi_1(M)$ contains no nontrivial normal abelian subgroup, and such that
		$$\injrad(M)\geq \varepsilon, \hspace{0.5 cm} \diam(M)\leq D.$$
		 Then for any group $\Gamma$ acting properly discontinuously by isometries on $X$, we have $\vol(X\slash\Gamma)\geq \mu$.
	\label{cor:minorbgen}
	\end{cornr}
	\begin{rmknr}\mbox{}
		\begin{enumerate}[(1)]
			\item  Note that if $X$ is a symmetric space then $X$ has a compact manifold quotient $M$ such that $\pi_1(M)$ contains no nontrivial normal abelian subgroups. Therefore Corollary \ref{cor:minorbgen} generalizes the theorem of Kazhdan-Margulis. However, as the proof uses the result of Kazhdan-Margulis, this does not provide a new proof.
 			\item The existence of minimal orbifolds is entwined with Conjecture \ref{conj:magicnr}. Namely, if the conjecture is true for a manifold $M$, then the universal cover has the minimal orbifolds property.
	 		\item Minimal \emph{manifolds} are related to collapsing in Riemannian geometry. If $M_k$ is a sequence of Riemannian $n$-manifolds with $\vol(M_k)\rightarrow 0$ but with uniformly bounded sectional curvatures and diameter, then $(M_k)_k$ is a collapsing sequence. By the work of Cheeger, Fukaya, and Gromov, such a manifold admits a nilpotent Killing structure along which the collapse occurs. This in turn forces topological restrictions. See \cite{chgrcoll1, chgrcoll2, chfugr, chrocoll1, chrocoll2} for more information.
	 	\end{enumerate}
	 \end{rmknr}
\subsection*{Outline of proofs} We give a brief outline of the proof of the Main Theorems \ref{thm:cmptbddbound} and \ref{thm:coverbddbound}. Suppose that there exist closed Riemannian manifolds $(M_k, g_k)$ satisfying the bounds of Equation \ref{eq:bounds} and with isometry groups $I_k$ such that $|I_k|\rightarrow\infty$. By a result of Anderson, the bounds of Equation \ref{eq:bounds} imply that along a subsequence $M_k$ are diffeomorphic to some closed manifold $M$ and $g_k\rightarrow g$ for some Riemannian metric $g$.  

Write $I:=\isom(M,g)$. Then $I$ is a Lie group acting by $C^2$ diffeomorphisms, so for Theorem \ref{thm:cmptbddbound} it suffices to show $I$ is infinite. We produce many elements of $I$ in the following way. Since $g_k\rightarrow g$, an isometry of $g_k$ is nearly an isometry of $g$. This suffices to show the family $\cup_k I_k$ is uniformly equicontinuous. Therefore any infinite sequence of elements of $\cup_k I_k$ subconverges to an element of $I$. The difficulty is to produce infinitely many \emph{distinct} such limits.

Let us just sketch how to produce one nontrivial element this way. The key tool here is a theorem of Newman, which shows the maximal diameter (with respect to a fixed metric) of all $I_k$-orbits is bounded away from 0 independently of $k$. So we can choose $\delta>0$ and $f_k\in I_k$ and $p_k\in M$ such that
	$$d(f_k p_k, p_k)\geq \delta.$$
This inequality passes to limits, so that the limit $f$ (along a subsequence) will be nontrivial.

For the proof of Theorem \ref{thm:coverbddbound}, there is the additional difficulty that we may have $p_k\rightarrow\infty$, so that no information about the limit is obtained. However, we show that we can choose $p_k$ in a compact subset of $\widetilde{M}$.

\subsection*{Outline of the paper} In Section \ref{sec:prelim} we discuss some preliminary tools that will be used in the proofs. In the next two sections we prove the main theorems. In Section \ref{sec:cmptbddbound} we first prove Theorem \ref{thm:coverbddbound}. Then we prove Theorems \ref{thm:cmptbddbound} and \ref{thm:quantborel}. In Section \ref{sec:coverreg} we prove Theorem \ref{thm:coverreg} that establishes a smooth symmetry gap for certain aspherical manifolds. We prove Corollary \ref{cor:minorbgen} on the existence of minimal orbifolds in Section \ref{sec:minorb}. Finally in Section \ref{sec:nordetect} we prove a technical result (Theorem \ref{thm:nordetect}) needed in the proof of Theorem \ref{thm:coverreg}.

\subsection*{Acknowledgments:} I am pleased to thank Max Engelstein and Katie Mann for helpful conversations regarding Section \ref{sec:isomlie}. I am grateful to Tu Tam Nguyen Phan for helpful discussions about piecewise locally symmetric spaces. Many thanks to Shmuel Weinberger for generous advice and countless helpful suggestions. I am very grateful to my thesis advisor Benson Farb for his invaluable advice and relentless enthusiasm during the completion of this work and extensive comments on an earlier version of this paper. I would like to thank the University of Chicago for support.

\section{Preliminaries}
\label{sec:prelim}

\subsection{Convergence of Riemannian manifolds} The discussion in this section is based on \cite[Chapter 10]{petersen}. We first define the important notion of convergence of Riemannian manifolds.
	\begin{dfnnr} Let $r\geq 1$ and let $(M_k, p_k, g_k)$ and $(M,p,g)$ be pointed complete $C^r$-Riemannian manifolds. For $r>0$ we say that
		$$(M_k, p_k, g_k)\rightarrow (M, p, g) \hspace{0.2 cm} \textrm{in the } C^r\textrm{-topology}$$
	if for every $R>0$ we can find a domain $\Omega\supseteq B_M(p;R)$ and embeddings
		$$f_k:\Omega\hookrightarrow M_k$$
	such that 
	\begin{enumerate}
		\item $B_{M_k}(p_k;R)\subseteq f_k(\Omega)$,
		\item $f_k^{\ast}g_k\rightarrow g$ in the $C^r$-topology on metrics on $\Omega$, and
		\item $f_k(p)=p_k$.
	\end{enumerate}
	\label{dfn:conv}\end{dfnnr}
	\begin{rmknr} In the above definition we allow $M_k$ and $M$ to be noncompact. If $M_k$ have universally bounded diameter, then it is easy to see that for $k\gg 1$ all maps $f_k$ are diffeomorphisms and basepoints can be chosen such that Condition (3) holds.\end{rmknr}
	\begin{rmknr} We will be especially interested in the case that $M_k$ is diffeomorphic to $M$ for all $k$, but the metrics $g_k$ are distinct. In this case, it is important to note that even if $(M, g_k)\rightarrow (M,g)$ in the $C^r$-topology, the metrics $g_k$ may not converge to $g$. For more information we refer to \cite{petersen}.\end{rmknr}	
There is a large amount of work on compactness results of families of Riemannian manifolds with certain geometric restrictions. This started with the result of Cheeger that the family of Riemannian manifolds with uniformly bounded sectional curvature, injectivity radius and diameter is precompact in the $C^{1,\alpha}$-topology \cite{cheegerfinite}. This was subsequently improved by Anderson to the following theorem that uses Ricci curvature instead of sectional curvature.
	\begin{thmnr}[Anderson \cite{metcmpt}] Let $\Lambda,\varepsilon,D>0$ and $n\geq 1$. Also fix $0<\alpha<1$. The family of closed Riemannian $n$-manifolds $(M,g)$ such that
		$$|\Ric(M_k, g_k)|\leq\Lambda, \hspace{0.5 cm} \injrad(M_k, g_k)\geq \varepsilon, \hspace{0.5 cm} \diam(M_k, g_k)\leq D$$
	is precompact in the $C^{1,\alpha}$-topology. In particular this family contains only finitely many diffeomorphism types.
	\label{thm:metcmpt}
	\end{thmnr}
	
\subsection{Isometry group of limit metric} \label{sec:isomlie} Motivated by Theorem \ref{thm:metcmpt}, we study the isometry groups of $C^{1,\alpha}$-metrics on Riemannian manifolds. Let $M$ be a connected smooth manifold (not necessarily compact) and let $g$ be a $C^r (r\geq 1)$ Riemannian metric on $M$. Here we mean that $r=k+\alpha$ for an integer $k\geq 1$ and $0<\alpha<1$ and when we write $g=g_{ij} dx^i dx^j$ in smooth local coordinates on $M$, we have that for every $i,j$, the $k$th derivatives of $g_{ij}$ are H\"older continuous with exponent $\alpha$. Set $I:=\isom(M,g)$. We will need the following result, which is probably well-known to the experts, but we could not find the proof in the literature for metrics of regularity less than $C^2$, so we will include it here.
\begin{propnr} The topology of uniform convergence on compact subsets of $M$ induces the structure of a Lie group on $I$ (possibly with infinitely many components), and the action of $I$ on $M$ is $C^{k+1}$.\label{prop:isomlie}\end{propnr}
\begin{proof} The Arzel\`a-Ascoli theorem implies that $I$ is locally compact. A theorem of Calabi-Hartman \cite{isomreg} implies that $I$ consists of $C^{r+1}$-diffeomorphisms of $M$. One should be warned here that the result in \cite{isomreg} is stated for all $r>0$, but there is a flaw in the proof for $r<1$ as was shown by Lytchak-Yaman \cite{isomnotreg}. However, the conclusion is still correct, even for $r<1$, as was shown by Taylor \cite{isomregafterall}.

Since $r\geq 1$, it follows that $I$ consists of $C^2$ diffeomorphisms of $M$. Further note that Bochner-Montgomery showed \cite[Theorem 1]{loccmptdiff} that a locally compact subgroup of Diff$^1(M)$ has the no small subgroups property. By the solution to Hilbert's fifth problem by Gleason, Montgomery-Zippin, and Yamabe (see \cite{hilbertsfifth}), it follows that $I$ is a Lie group. 

Finally, Montgomery showed \cite[Corollary 1]{actionreg} that whenever a Lie group acts effectively by $C^k$ diffeomorphisms on a manifold, then the action is $C^k$.\end{proof}
\subsection{Equivariant Gromov-Hausdorff convergence} \label{sec:egh} For information about Gromov-Hausdorff convergence we refer to \cite{riemnonriem}, \cite{bubuiv}, and \cite{fukayasurvey}. An equivariant version of Gromov-Hausdorff convergence was developed by Fukaya \cite{fukayaegh}. We introduce the following notation. If $(X, p)$ is a pointed proper metric space and $G$ is a closed subgroup of isometries of $X$, we set for $R>0$
	$$G(R):=\{g\in G\mid d(p, gp)<R\}.$$
Now we make the following definition.
	\begin{dfnnr} Let $(X, p), (Y,q)$ be pointed proper metric spaces, and let $G\subseteq \isom(X)$ and $H\subseteq \isom(Y)$. For $\varepsilon>0$, an $\varepsilon$-\emph{equivariant Gromov Hausdorff approximation} is a triple $(f,\varphi,\psi)$ where
		\begin{enumerate}
			\item $f:B_X(p;\frac{1}{\varepsilon})\rightarrow Y$,
			\item $\varphi: G(\frac{1}{\varepsilon})\rightarrow H(\frac{1}{\varepsilon})$, and
			\item $\psi: H(\frac{1}{\varepsilon})\rightarrow G(\frac{1}{\varepsilon}),$
		\end{enumerate}
	satisfying
		\begin{enumerate}
			\item $f(p)=q$,
			\item $B_Y(q;\frac{1}{\varepsilon})\subseteq N_\varepsilon(f(B_X(p;\frac{1}{\varepsilon})))$ where $N_\varepsilon$ denotes the $\varepsilon$-neighborhood,
			\item for $x,y\in B_X(p;\frac{1}{\varepsilon})$, we have
				$$|d_Y(f(x), f(y))-d_X(x,y)|<\varepsilon,$$
			\item for $g\in G(\frac{1}{\varepsilon})$ and $x\in B_X(p;\frac{1}{\varepsilon})$ such that $gx\in B_X(p;\frac{1}{\varepsilon})$, we have
				$$d_Y(f(gx), \varphi(g)f(x))<\varepsilon,$$
			\item for $h\in H(\frac{1}{\varepsilon})$ and $x\in B_X(p;\frac{1}{\varepsilon})$ such that $\psi(h)x\in B_X(p;\frac{1}{\varepsilon})$, we have
				$$d_Y(f(\psi(h)x), h f(x))<\varepsilon.$$
		\end{enumerate}
	The \emph{equivariant Gromov-Hausdorff distance} $d_{\textrm{eGH}}((X,G,p),(Y,H,q))$ is the infimum of $\varepsilon$ such that there exists an $\varepsilon$-equivariant Gromov-Hausdorff approximation from $(X,G,p)$ to $(Y, H, q)$ and vice versa. This defines an obvious notion of convergence.
	\end{dfnnr}
	\begin{rmknr} It is \emph{not} required that $f$ is continuous, or that $\varphi, \psi$ are (restrictions of) homomorphisms.\label{rmk:nohomegh}\end{rmknr}
The following result by Fukaya-Yamaguchi relates ordinary Gromov-Hausdorff convergence to the equivariant case.
	\begin{thmnr}[{\cite[Prop 3.6]{fyalmostneg}}] Let $(X_k, p_k)$, $(X, p)$ such that $(X_k, p_k)\rightarrow (X,p)$and let $G_k\subseteq \isom(X_k)$ be closed subgroups. Then there exists a closed subgroup $G\subseteq \isom(X)$ such that $(X_k, G_k, p_k)$ subconverges in the equivariant Gromov-Hausdorff topology to $(X, G, p)$.\label{thm:eghconv}\end{thmnr}
We also have the following relation between equivariant Gromov-Hausdorff convergence and Gromov-Hausdorff convergence of the orbit spaces.
	\begin{thmnr}[{\cite[Theorem 2.1]{fukayaegh}}] Let $(X_k, p_k)$ and $(X,p)$ be pointed proper metric spaces and $G_k\subseteq \isom(X_k)$ and $G\subseteq \isom(X)$ closed subgroups such that $(X_k, G_k, p_k)\to(X, G, p)$ in the equivariant Gromov-Hausdorff sense. Then 
		$$(X\slash G_k, [p_k])\to(X\slash G, [p]).$$
	\label{thm:orbitgh}
	\end{thmnr}
However, if $(X_k, G_k, p_k)\to(X, G, p)$, it is in general very difficult to relate the group structure of $G_k$ to that of $G$, essentially because equivariant Gromov-Hausdorff distance makes no reference to morphisms. Still we have the following result of Fukaya-Yamaguchi that detects suitable normal subgroups in the limit along the sequence.
	\begin{thmnr}[{\cite[Theorem 3.10]{fyalmostneg}}] Let $(X_k, p_k)$ and $(Y,q)$ be pointed proper metric spaces. Let $\Gamma_k\subseteq\isom(X_k)$ and $G\subseteq\isom(Y)$ be closed subgroups such that 
		$$(X,\Gamma_k,p_k)\to(Y, G, q).$$
	Let $G'\subseteq G$ be a normal subgroup and assume the following.
		\begin{enumerate}
			\item $G\slash G'$ is discrete and finitely presented,
			\item $Y\slash G$ is compact,
			\item $\Gamma_k$ acts on $X_k$ properly discontinuously and freely,
			\item $X_k$ is simply-connected and there exists $R_0$ such that $\pi_1(B_Y(q;R_0))\rightarrow \pi_1(Y)$ is surjective, and
			\item there exists $R_1$ such that $G'(R_1)$ generates $G'$.
		\end{enumerate}
	Then there exist normal subgroups $\Gamma_k'\subseteq\Gamma_k$ such that
		\begin{enumerate}
			\item $(X_k, \Gamma_k', p_k)\to(Y, G', q)$,
			\item $\Gamma_k\slash\Gamma_k'$ and $G\slash G'$ are isomorphic for $k$ sufficiently large, and
			\item there exists $R_2$ such that for all $k$ sufficiently large, $\Gamma_k'$ is generated by $\Gamma_k'(R_2)$.
		\end{enumerate}
	\label{thm:fynordetect}
	\end{thmnr}
Fukaya-Yamaguchi remarked that the conclusions of Theorem \ref{thm:fynordetect} might remain true without the assumption that $\Gamma_k$ act freely on $X_k$. We obtain the following towards this generalization.
	\begin{thmnr} Assume the hypotheses of Theorem \ref{thm:fynordetect} hold except $\Gamma_k$ are not assumed to act freely on $X_k$. Assume in addition that $X_k$ is a manifold for every $k$. Then the conclusion of Theorem \ref{thm:fynordetect} still holds.\label{thm:nordetect}\end{thmnr}
We will prove Theorem \ref{thm:nordetect} in Section \ref{sec:nordetect}. However, let us briefly explain the idea. Fukaya-Yamaguchi used the assumption of free actions to define certain covering spaces, and $\Gamma_k'$ will be the fundamental group of such a covering space. If the action of $\Gamma_k$ on $X_k$ is not free, then this construction will not yield covering spaces, and we cannot define $\Gamma_k'$ in this way. However, if $X_k$ are assumed to be manifolds then the spaces of Fukaya-Yamaguchi are naturally orbifolds. In fact we can show they are good (sometimes also called `developable') orbifolds. For good orbifolds, Thurston developed a theory of orbifold covering spaces. We can then define $\Gamma_k'$ to be the orbifold fundamental group of a suitable orbifold covering space, and the rest of the proof of Fukaya-Yamaguchi carries through verbatim.

\subsection{Farb-Weinberger's work on symmetries of universal covers} \label{sec:fw} We review some of the material of \cite{FW} since we will use some of the results and variations on their arguments. The main theorem is as follows.
	\begin{thmnr}[Farb-Weinberger] Let $M$ be a closed, aspherical Riemannian manifold. Then either $[\isom(\widetilde{M}):\pi_1(M)]<\infty$, or $M$ is isometric to a Riemannian orbibundle
		$$F\rightarrow M\rightarrow B$$
	where
		\begin{itemize}
			\item $B$ is a good Riemannian orbifold, and
			\item each fiber $F$ is isometric (with respect to the induced metric) to a nontrivial closed, aspherical locally homogeneous space.
		\end{itemize}
	\label{thm:fw}
	\end{thmnr}
Here a \emph{Riemannian orbibundle} is a map that is locally modelled on the quotient map $V\times_G F\rightarrow V\slash G$ for a finite group $G$ acting on a fixed smooth manifold $F$ and $(V,G)$ is a chart for the orbifold $B$. Further we require that $G$ acts isometrically on $V\times F$ and the projection to $V$ is a Riemannian submersion. We will need the following two useful facts obtained in the course of the proof. 

We will fix the following notation for the rest of the section. Let $M$ be a closed Riemannian manifold with a $C^1$-Riemannian metric. Set $I:=\isom(\widetilde{M})$ and $\Gamma:=\pi_1(M)$. Assume that $[I:\Gamma]=\infty$. By Proposition \ref{prop:isomlie} we know that $I$ is a Lie group (possibly with infinitely many components) and $I$ acts on $\widetilde{M}$ by $C^2$ diffeomorphisms. Set $\Gamma_0:=\Gamma\cap I^0$. Then we have the following.
\begin{lemnr}[{\cite[Claims I and II]{FW}}] \mbox{}
	\begin{enumerate}[(1)]
		\item $\Gamma_0$ is a cocompact lattice in $I^0$.\label{lem:lattice}
		\item If $M$ is aspherical then $I^0$ contains no nontrivial compact factors.\label{lem:nocmpt}
	\end{enumerate}
	\label{lem:fwclaims}
	\end{lemnr}
\begin{rmknr} Farb-Weinberger prove Lemma \ref{lem:fwclaims} only for smooth Riemannian manifolds, but their proof works verbatim for $C^1$ metrics as well.\end{rmknr}
Assume in addition that $M$ is aspherical and $\pi_1(M)$ contains no normal abelian subgroup. Then we get far stronger results. First we see that the structure of $I^0$ is very constrained as follows.
\begin{propnr}[{\cite[Proposition 3.3]{FW}}] Let $M$ be aspherical such that $\pi_1(M)$ contains no nontrivial normal abelian subgroup. Then $I^0$ is semisimple with finite center and no compact factors.\label{prop:ss}\end{propnr}
Now assume in addition that $(M,g)$ is smooth. Then Proposition \ref{prop:ss} implies that the fiber $F$ obtained in Theorem \ref{thm:fw} is a nonpositively curved locally symmetric space. Using the theory of harmonic maps for nonpositively curved manifolds, we can construct a section of the orbibundle obtained in Theorem \ref{thm:fw}. Together this yields the following result.
\begin{thmnr}[{\cite[Proposition 3.1]{FW}}] Let $M$ be aspherical such that $\pi_1(M)$ contains no nontrivial normal abelain subgroup. Then a finite cover of $M$ is a Riemannian warped product $B\times_f N$ where $N$ is a locally symmetric space of noncompact type and $f:B\rightarrow \bbR_{>0}$ is a smooth function.\label{thm:warped}\end{thmnr}
Recall that if $(X, g_X)$ and $(Y, g_Y)$ are Riemannian manifolds, then the Riemannian warped product $X\times_f Y$ has underlying manifold $X\times Y$ equipped with the metric
	$$g|_{(x,y)}=g_X|_x \oplus f(x) g_Y|_y.$$
The first step in the proof of Theorem \ref{thm:warped} is the following proposition. Since the proof will be useful for us, we give a sketch.
\begin{propnr} There is a finite index subgroup $\Gamma'\subseteq\Gamma$ such that 
	\begin{itemize}
		\item $\Gamma_0\subseteq \Gamma'$, and
		\item $\Gamma'\cong \Gamma_0\times (\Gamma'\slash\Gamma_0)$.
	\end{itemize}
	\label{prop:prod}
\end{propnr}
\begin{proof} Consider the short exact sequence
	$$1\rightarrow I^0\rightarrow \langle \Gamma, I^0\rangle \rightarrow \Gamma\slash\Gamma_0\rightarrow 1.$$
This short exact sequence gives rise to a morphism 
	$$\rho:\Gamma\slash\Gamma_0\rightarrow\Out(\Gamma_0).$$
Since Out$(I^0)$ is finite we can assume $\rho$ is trivial by passing to a finite index subgroup of $\Gamma$. Now consider the extension
	\begin{equation}1\rightarrow \Gamma_0\rightarrow \Gamma\rightarrow \Gamma\slash\Gamma_0\rightarrow 1.\label{eq:fwses}\end{equation}
This extension gives rise to a morphism
	$$\sigma: \Gamma\slash\Gamma_0\rightarrow \Out(\Gamma_0).$$
Let $N_{I^0}\Gamma_0$ be the normalizer of $\Gamma_0$ in $I^0$. Since we know $\Gamma\slash\Gamma_0$ acts by inner automorphisms on $I^0$, it follows that $\sigma$ has image in the finite group $N_{I^0}\Gamma_0\slash\Gamma_0$. Therefore by passing to a finite index subgroup of $\Gamma$, we can assume that $\sigma$ is trivial. The extension \ref{eq:fwses} is now determined by a cohomology class in $H^2(\Gamma\slash\Gamma_0, Z(\Gamma_0))$. But since $I^0$ is a connected semisimple Lie group with finite center and $\Gamma_0$ is a torsion-free cocompact lattice in $I^0$, we know that $Z(\Gamma_0)=1$. Therefore the extension \ref{eq:fwses} is trivial, so that
	$$\Gamma\cong\Gamma_0\times (\Gamma\slash\Gamma_0).$$
\end{proof}

\section{Isometries of Riemannian manifolds and covers}
\label{sec:cmptbddbound}
The goal of this section is to prove our main theorems. First we prove Theorem \ref{thm:cmptbddbound} and \ref{thm:coverbddbound}. In fact both follow from a general theorem about isometries of covers (Theorem \ref{thm:gencover} below) which we prove first. Theorem \ref{thm:cmptbddbound} (resp. \ref{thm:coverbddbound}) follows from the case of the trivial cover $M\rightarrow M$ (resp. the universal cover $\widetilde{M}\rightarrow M$) combined with the Anderson Metric Compactness theorem (Theorem \ref{thm:metcmpt}). At the end of the section we prove Theorem \ref{thm:quantborel}. 
	\begin{thmnr} Let $n\geq 1$ and $\Lambda, D, \varepsilon>0$. Let $M$ be a closed smooth manifold and $M'\rightarrow M$ any regular cover of $M$, and let $\Gamma$ denote the deck group. Let $g_k$ be Riemannian metrics on $M$ such that
			\begin{equation}
				|\Ric_{g_k}|\leq \Lambda, \hspace{0.5 cm} \injrad(M,g_k)\geq \varepsilon, \hspace{0.5 cm} \diam(M,g_k)\leq D.
			\end{equation}
	and $[\isom(M', g_k):\Gamma]\rightarrow\infty$. Then $M'$ admits an effective proper $C^2$ action by a positive-dimensional Lie group $G$ containing $\Gamma$.
	\label{thm:gencover}
	\end{thmnr}
\begin{proof}[Proof of Theorem \ref{thm:gencover}] By Anderson's Metric Compactness Theorem (Theorem \ref{thm:metcmpt}) we can assume there exist diffeomorphisms $f_k:M'\rightarrow M'$ and a $C^{1,\alpha}$ Riemannian metric $g$ on $M'$ (for some fixed $0<\alpha<1$) such that $f_k^\ast g_k\rightarrow g$ in the $C^{1,\alpha}$-topology. Since
	$$I(M', g_k) \cong I(M', f_k^\ast g_k)$$
we can assume that $f_k=\id$ for all $k$, so that $g_k\rightarrow g$ in the $C^{1,\alpha}$-topology.

Write $I_k:=\isom(M',g_k)$ for the isometry groups of $M'$ and $I:=\isom(M',g)$ for the isometry group of the limit metric. By Proposition \ref{prop:isomlie}, we have that $I$ is a Lie group (possibly with infinitely many components) acting by $C^{2,\alpha}$ diffeomorphisms of $M'$. 

Therefore it suffices to show that $[I:\Gamma]=\infty$. We prove by contradiction that $I\slash\Gamma$ contains infinitely many cosets $h\Gamma$ such that $h=\displaystyle\lim_{k\rightarrow\infty} h_k$ where $h_k\in I_k$. The following lemma allows us to produce many isometries.
	\begin{lemnr}[Isometry precompactness] $\mathcal{I}:=\cup_k I_k$ is uniformly equicontinuous on compact subsets of $M'$. Hence if $f_k\in I_k$ such that for some $p\in M'$ the sequence $\{f_k(p)\}_k$ is bounded, then $(f_k)_k$ subconverges to some $f\in I$. \label{lem:eqcnts}\end{lemnr}
	\begin{proof} The second claim follows from the first claim combined with the Arzel\`a-Ascoli theorem. The fact that the limit is an isometry of $g$ follows easily from $g_k\rightarrow g$. We prove the first claim. Let $K\subseteq M'$ be any compact subset and $\varepsilon>0$. Then for $k$ sufficiently large, we have $|d(x,y)-d_k(x,y)|<\frac{\varepsilon}{3}$ for any $x,y\in K$. Fix such $k$ and choose $x,y\in K$ with $d(x,y)<\frac{\varepsilon}{3}$. Then for any $f\in I_k$, we have
			\begin{align*}
			d(f(x),f(y))&\leq d_k(f(x),f(y))+\frac{\varepsilon}{3}\\
						&=d_k(x,y)+\frac{\varepsilon}{3}\\
						&\leq d(x,y)+\frac{2\varepsilon}{3}\\
						&<\varepsilon.
			\end{align*}		\end{proof}
Suppose now there are only finitely many cosets in $I\slash\Gamma$ represented by elements that are limits $\displaystyle\lim_{k\rightarrow\infty} h_k$ with $h_k\in I_k$. Choose such representatives $h^0, \dots, h^r$ with $h^j=\displaystyle\lim_{k\rightarrow \infty} h^j_k$ in the Gromov-Hausdorff topology where $h^j_k\in I_k$ and $h^0=e$. We can assume that $H:=\cup_j h^j\Gamma$ is a subgroup of $I$. For $k\geq 1$, set 
	$$H_k:=\bigcup_{0\leq j\leq r} h^j_k \Gamma.$$
Note that we do not assume that $H_k$ is a subgroup of $I_k$. Our goal is to produce an isometry $f\in I\backslash H$.

Since $\textrm{Fix}(h)$ is nowhere dense for every $h\neq e$, we can choose $p\in M'$ not fixed by any $h\in H\backslash\{e\}$. Since $H$ contains $\Gamma$ as a finite index subgroup, $H$ acts properly discontinuously on $M'$. Therefore we can choose $\eta>0$ such that 
	$$\eta<\frac{1}{4}\min_{e\neq h\in H} d(p, h (p)).$$
In addition choose $\eta<\frac{\varepsilon}{4}$, where $\varepsilon$ is the lower bound on the injectivity radii of $(M, g_k)$. We consider two cases, depending on whether the $I_{k_l}$-orbit of $p$ lies in an $\eta$-neighborhood of the $H_{k_l}$-orbit of $p$ or not. We refer to these as the case of a `concentrated orbit' and `diffuse orbit'.
\subsubsection*{Case 1 (diffuse orbit)} Assume that for a subsequence $k_l\rightarrow \infty$, we have
	$$I_{k_l}\cdot p\nsubseteq B_{k_l}(H_{k_l} \cdot p; \eta).$$
Then we can choose $f_{k_l}\in I_{k_l}$ such that $d_{k_l}(f_{k_l}(p), h (p))\geq \eta$ for all $h\in H_k$. By postcomposing $f_{k_l}$ by an element of $\Gamma$ we can assume that 
	$$d_{k_l}(f_{k_l} (p), p)\leq \diam(M,g_{k_l})\leq D.$$
Therefore along a subsequence we have $f_{k_l}\rightarrow f$ for some $f\in I$, and we know
	$$d(f(p), h (p))=\lim_{l\rightarrow \infty} d_{k_l}(f_{k_l}(p),h_{k_l} (p))\geq \eta$$
for every $h\in H$. It follows that $f\neq h$ for any $h\in H$, as desired.
\subsubsection*{Case 2 (concentrated orbit)} For $k\gg 1$ we have
	$$I_k \cdot p\subseteq \bigcup_{h\in H} B_k(h (p);\eta).$$
Set
	$$\Lambda_k:=\{f\in I_k\mid d_k(f(p), p)<\eta\}.$$
\textbf{Step 1 ($\Lambda_k$ is a group for $k\gg 1$)}. We claim that $\Lambda_k$ is a finite subgroup of $I_k$ for $k\gg 1$. It is clear that whenever $f\in \Lambda_k$, we have $f^{-1}\in \Lambda_k$ and that $\Lambda$ contains $\id_{M'}$. It remains to show that $\Lambda_k$ is closed under multiplication. Let $f_1, f_2\in \Lambda_k$. Then there exists $0\leq j\leq r$ such that
		\begin{align} d_k(f_1 f_2 (p), h^k_j (p))<\eta. \label{eq:smalldisp} \end{align}
	Our goal is to show $j=0$. Note that if $1\leq i\leq r$ then
		$$\lim_{k\rightarrow \infty} d_k(h^k_i (p), p) = d(h_i (p), p)\geq 4\eta. $$
	So we can choose $k\gg 1$ such that for any $1\leq i\leq r$, we have
		\begin{align}d_k(h^i_k (p), p)\geq 3\eta.\label{eq:bigdisp} \end{align}
	It follows from Equations \ref{eq:smalldisp} and \ref{eq:bigdisp} that for any $j\neq 0$, we have
		\begin{align} d_k(f_1 f_2 (p), p)\geq 2\eta. \label{eq:prop1} \end{align}
	On the other hand, since $f_1^{-1}, f_2\in\Lambda_k$, we find
		\begin{align*}
			d_k(f_1 f_2 (p), p)&=d_k (f_2 (p), f_1^{-1} (p))\\
					 &\leq d_k(f_2 (p), p) + d_k(p, f_1^{-1}(p))\\
					 &< 2\eta,
		\end{align*}
	contradicting Equation \ref{eq:prop1}. Therefore we must have $j=0$. It remains to show that $\Lambda_k$ is finite. This follows immediately from the definition of $\Lambda_k$ and the proper discontinuity of the action of $I_k$ on $M'$. This completes Step 1.
	
\textbf{Step 2 (Nontriviality)}. We show that $\Lambda_k$ is nontrivial for $k\gg 1$. Suppose $\Lambda_{k_l}$ is trivial for some subsequence $k_l\rightarrow\infty$. Let $f\in I_{k_l}$. By assumption there exists $0\leq j\leq r$ and $h\in h^j_{k_l}\Gamma$ such that $d(f(p), h(p))<\eta$. It follows that $f^{-1} h\in\Lambda_{k_l}=\{e\}$, so $f=h\in h^j_{k_l}\Gamma$ Since $f\in I_{k_l}$ was arbitrary, we conclude that $[I_{k_l}:\Gamma]\leq r+1$. This is a contradiction for $l\gg 1$, which completes Step 2.

Newman proved that there exists $\delta>0$ such that such for any nontrivial compact group $G$ acting effectively on $M'$, there exists a $G$-orbit with diameter at least $\delta$ with respect to the metric $d$ \cite[Corollary III.9.6]{brtrgps}. In particular there exists $\delta>0$ (independent of $k$) and $q_k\in M'$ such that $\diam(\Lambda_k \cdot q_k)\geq \delta$. A major difficulty is that we could have $q_k\rightarrow\infty$, so that the inequality $d(f_k (q_k), q_k)\geq \delta$ does not give any information as $k\rightarrow\infty$. We resolve this in the following way.

\textbf{Step 3 (Trapping $q_k$)} We prove there exists a compact subset $K$ of $M'$ such that for every $k\gg 1$, there exists $e\neq f\in\Lambda_k$ with a fixed point in $K$.

Let $k\gg 1$ such that Step 1 holds. Consider
	$$V_k:=\bigcup_{\lambda\in\Lambda_k} \overline{B_k(\lambda (p);\eta)}.$$
Note that 
	$$\diam_k V_k\leq 3\eta\leq \textrm{injrad}(M,g_k).$$
Therefore $U_k$ is contained in a geodesically convex ball $B_k$ centered at $p$, so that every pair of points of $V_k$ is joined by a unique geodesic with respect to $g_k$. Let $C_k$ be the convex hull of $U_k$. Since $V_k$ is a $\Lambda_k$-invariant set and $\Lambda_k$ acts by isometries with respect to $g_k$, it follows that $C_k$ is also $\Lambda_k$-invariant. Further, since $V_k$ is closed and convex, it follows that $C_k$ is homeomorphic to an $n$-dimensional disk.

Now let $e\neq f\in \Lambda_k$. By possibly replacing $f$ by a power, we can assume $f$ has prime order $\ell$. Consider the action of $\langle f\rangle$ on $C_k$. Since $C_k$ is contractible, it is well-known that $f$ must have a fixed point on $C_k$, for otherwise $C_k\slash\langle f\rangle$ would be a finite dimensional $K(\bbZ\slash \ell\bbZ, 1)$, which is impossible. Let $x_k\in C_k$ be fixed by $f$.

Now set $K:=\overline{B(p;4\eta)}$. Note that by construction, we have $d_k(x_k,p)\leq 3\eta$. Therefore for $k\gg 1$, we have $x_k\in K$, which proves Step 2.

\textbf{Step 4 (Constructing $q_k$)}. Choose $e\neq f_k\in\Lambda_k$ with fixed points $x_k\in K$ as above. Since $K$ is compact, after passing to a subsequence, we can assume that $x_k\rightarrow x$ for some $x\in M'$. Now consider the action of $\langle f_k\rangle$ on $B_k(x_k;\eta)$. 

%

Let $\nabla_{x_k, x}^{g_k}$ denote parallel transport from $x_k$ to $x$ with respect to the metric $g_k$. Then $\nabla_{x_k,x}^{g_k}\circ \left(\exp^{g_k}_{x_k}\right)^{-1}$ is a diffeomorphism from $B_k(x_k;\eta)$ to the $\eta$-ball $B_k:=B_k(0;\eta)\subseteq T_x M'$ with respect to $g_k$. The above diffeomorphism defines a smooth conjugacy between the action of $\langle f_k\rangle$ on $B_k(x_k;\eta)$ and an action on $B_k$.

Choose linear isomorphisms $A_k$ of $T_x M'$ such that $A_k$ conjugates the quadratic form $g_k|_x$ to $g|_x$. Since $g_k\rightarrow g$ in the $C^{1,\alpha}$-topology, we have $g_k|_x\rightarrow g|_x$, so that we can choose $A_k$ such that $A_k\rightarrow \id$. Denote by $B:=B(0;\eta)$ the ball in $T_xM'$ with respect to the metric $g$. Then $A_k$ conjugates the action of $\langle f_k\rangle$ on $B_k$ to an action on $B$.

Now let $\delta>0$ be the constant from Newman's theorem for the metric space $(B,d)$. By choice of $\delta$, there exists $v_k\in B$ and $l_k\geq 1$ such that
	$$d(v_k, f^{l_k}_k (v_k))\geq \delta.$$
Along a subsequence, we can assume $v_k\rightarrow v$ for some $v\in \overline{B}$ and $f^{l_k}_k\rightarrow f$. Write
	$$q_k:=\exp^{g_k}_{x_k} \nabla^{g_k}_{x, x_k} A_k^{-1} v_k$$
where the exponential map and parallel transport are with respect to $g_k$. Similarly define $q:=\exp^g_x v$ where the exponential map is with respect to the metric $g$.

\textbf{Step 4 (Convergence)} We claim that $q_k\rightarrow q$. First note that
		$$d(q_k, \exp_x^{g_k} A_k^{-1}v_k)\rightarrow 0$$
	since $x_k\rightarrow x$ (so that $\nabla^{g_k}_{x, x_k}\rightarrow \id$). Secondly,
		$$d(\exp^{g_k}_x A_k^{-1} v_k, \exp^g_x A_k^{-1} v_k)\rightarrow 0$$
	since $g_k\rightarrow g$ (so that $\exp^{g_k}_x \rightarrow \exp^g_x$). Combining these observations, we see that $q_k\rightarrow q$ if and only if $A_k^{-1} v_k\rightarrow v$. But the latter is obvious since $v_k\rightarrow v$ and $A_k\rightarrow\id$.

\textbf{Step 5 (Constructing a new isometry)}. We produce an isometry $f\in I$ such that $f\notin H$. Let $k\gg 1$ such that Steps 1, 2, and 3 hold.  Recall that $q_k\in M'$ have the property $\diam(\Lambda_k \cdot q_k)\geq \delta$ (where $\delta$ does not depend on $k$). Choose $f_k\in\Lambda_k$ such that $d(f_k (q_k), q_k)\geq \delta$. Along a subsequence we can assume that $q_k\rightarrow q$ and $f_k\rightarrow f$ for some $q\in M$ and $f\in I$. We claim that $f\notin H$.

By Step 4, we have
	$$d(f(q) , q)=\lim_{k\rightarrow \infty} d_k(f_k (q_k), q_k),$$
and by the same arguments as in the proof of Step 4, we have
	$$\lim_{k\rightarrow\infty} d_k(f_k (q_k), q_k)= \lim_{k\rightarrow \infty} d(f_k(\exp^g_x(A_k^{-1} v_k)), \exp^g_x(A_k^{-1} v_k)).$$
(Here the exponential maps are with respect to $g$). Since $g$ is a $C^{1,\alpha}$-metric with $\alpha>0$, exponential maps of $g$ are H\"older continuous with exponent $\alpha$. Therefore there is a constant $L>0$ such that
	$$\lim_{k\rightarrow \infty} d(f_k (\exp^g_x(A_k^{-1} v_k)), \exp^g_x(A_k^{-1} v_k))\geq L \lim_{k\rightarrow\infty} \left(\|f_k (v_k) - v_k\|_{g_x}\right)^{\frac{1}{\alpha}}\geq L\delta^{\frac{1}{\alpha}}>0.$$
It follows that $d(f(q), q)\geq L\delta^{\frac{1}{\alpha}}>0$ so $f$ is nontrivial. Further observe that
	$$d(f (p), p)=\lim_{k\rightarrow\infty} d_k(f_k (p), p)\leq \eta$$
and for $e\neq h\in H$ we have
	$$d(h (p), p)\geq 4\eta.$$
This completes Step 5, which shows that $f\notin H$. On the other hand by construction $f$ is a limit of a subsequence of elements of $f_k\in I_k$. This contradicts the definition of $H$. This is the desired contradiction that completes the proof of Theorem \ref{thm:gencover}.\end{proof}
\begin{proof}[Proof of Theorem \ref{thm:cmptbddbound}] Let $n\geq 1$ and $\Lambda, D, \varepsilon>0$. Suppose that for $k\geq 1$ there exists a smooth $n$-manifolds $M_k$ and a Riemannian metric $g_k$ on $M$ such that
	$$|\Ric_{g_k}|\leq \Lambda, \hspace{0.5 cm} \diam(M_k,g_k)\leq D, \hspace{0.5 cm} \injrad(M_k,g_k)\geq \varepsilon$$
and $|\isom(M,g_k)|\rightarrow\infty$. We want to show that for any $k\gg 1$, the manifold $M_k$ admits a $C^2$ circle action. By Anderson's Metric Compactness Theorem (Theorem \ref{thm:metcmpt}) we can assume that there is a closed smooth manifold $M$ such that $M_k$ is diffeomorphic to $M$ for all $k$. Theorem \ref{thm:gencover} applied to the trivial cover $M\rightarrow M$ implies that $M$ admits a $C^2$ circle action.\end{proof}
\begin{proof}[Proof of Theorem \ref{thm:coverbddbound}] Let $n\geq 1$ and $\Lambda, D, \varepsilon>0$. Suppose that for $k\geq 1$ there exists a smooth $n$-manifolds $M_k$ and a Riemannian metric $g_k$ on $M$ such that
	$$|\Ric_{g_k}|\leq \Lambda, \hspace{0.5 cm} \diam(M_k,g_k)\leq D, \hspace{0.5 cm} \injrad(M_k,g_k)\geq \varepsilon$$
and $[\isom(\widetilde{M}_k,\widetilde{g}_k):\pi_1(M_k)]\rightarrow\infty$. We want to show that for any $k\gg 1$, there is an effective proper $C^2$ action by a positive-dimensional Lie group $G_k$ containing $\pi_1(M_k)$ on $\widetilde{M}_k$ . By Anderson's Metric Compactness Theorem (Theorem \ref{thm:metcmpt}) we can assume that there is a closed smooth manifold $M$ such that $M_k$ is diffeomorphic to $M$ for all $k$. Theorem \ref{thm:gencover} applied to the universal cover $\widetilde{M}\rightarrow M$ implies that $\widetilde{M}$ admits a proper $C^2$ action by a positive-dimensional Lie group.\end{proof}
\begin{proof}[Proof of Theorem \ref{thm:quantborel}] Let $(M,g)$ be a closed aspherical manifold with centerless fundamental group and satisfying the bounds of Equation \ref{eq:bounds}. It is well-known that the lower bound on Ricci curvature $\Ric\geq\Lambda$ and $\diam\leq D$ give a packing inequaltiy. Namely, the number $N$ of disjoint $\frac{\varepsilon}{4}$-balls in $M$ is bounded only depending on $\Lambda$ and $D$. More explicitly, we have the following Bishop-Gromov inequality for $R>r$ and $p\in \widetilde{M}$
	$$\frac{\vol\left(B(p;R\right)}{\vol(B(p;r))}\leq \frac{V_k(R)}{V_k(r)}$$
where $V_k(s)$ is the volume of a ball of radius $s$ in a simply-connected space of constant curvature $k=\frac{-\Lambda}{n-1}$ and dimension $n$. Using that $M$ has diameter $\leq D$, it is easily follows that the number of disjoint $\frac{\varepsilon}{4}$-balls in $M$ is at most $N=\frac{V_k\left(\frac{D}{2}\right)}{V_k\left(\frac{\varepsilon}{4}\right)}.$ Take a maximal collection of disjoint $\frac{\varepsilon}{4}$-balls. Then the collection of balls with the same centers but radius $\frac{\varepsilon}{2}$ covers $M$.

Now let $X$ be the nerve associated to this cover (i.e. the vertices of $X$ are the balls $B_i, 1\leq i\leq m$, and two vertices are joined if and only if $B_i\cap B_j\neq \emptyset$, and $B_i, B_j, B_k$ span a triangle if and only if $B_i\cap B_j\cap B_k\neq \emptyset$). 

Since injrad$(M,g)\geq \varepsilon$, we know that each of these balls are embedded and all possible intersections are contractible.  It follows that $\pi_1(M)\cong \pi_1(X)$. Now let $f:M\rightarrow M$ be an isometry. Then define a map $F:X\rightarrow X$ as follows. Let $x\in X^0$ and $B$ the corresponding ball with center $p\in M$. Let $q\in M$ be a point that is the center of some ball in $\mathcal{U}$ and minimizes distance to $f(p)$. Note that there might be many such points, but in which case we arbitrarily choose one of them. Set $F(x):=q$. Then extend $F$ to a simplicial map $X\rightarrow X$.

It is easy to see that $F_\ast=f_\ast$ on fundamental groups (after making the identification $\pi_1(X)\cong \pi_1(M)$ induced by the inclusion $X\hookrightarrow M$). However, there are only finitely many possibilities for $F$. Explicitly the number is bounded above by the number of maps $X^0\rightarrow X^0$, which is $N^N=C$. Therefore there are at most $C$ possibilities for $f_\ast:\pi_1(M)\rightarrow \pi_1(M)$ up to conjugation. 

The theorem of Borel on group actions closed aspherical manifolds with centerless fundamental groups mentioned in Section \ref{sec:intr}, has the following more precise formulation \cite{borelnoact}. Let $G$ be a compact group acting on a closed aspherical manifold $N$ with centerless fundamental group. Then the map
	$$G\rightarrow \textrm{Out}(\pi_1 N)$$
is an embedding. Setting $G:=\isom(M)$ and $N:=M$ we obtain $|\isom(M)|\leq C$.\end{proof}

\section{Isometry groups of contractible manifolds}
\label{sec:coverreg}
In this section we prove Theorem \ref{thm:coverreg}.
\begin{proof}[Proof of Theorem \ref{thm:coverreg}] Suppose there exist closed, aspherical, smoothly irreducible Riemannian $n$-manifolds $(M_k,g_k)$ such that $\pi_1(M)$ contains no nontrivial normal abelian subgroup, and such that
	$$|\Ric_{g_k}|\leq\Lambda, \hspace{0.5 cm} \injrad(M_k,g_k)\geq \varepsilon, \hspace{0.5 cm} \diam(M_k,g_k)\leq D,$$
and $[I(\widetilde{M_k},\widetilde{g_k}):\pi_1(M_k)]\rightarrow\infty$. As before we can assume that $M_k$ are diffeomorphic to a single manifold $M$ and $g_k\rightarrow g$ in the $C^{1,\alpha}$-topology. By the work of Farb-Weinberger, if $[I(\widetilde{M},\widetilde{g}_k):\pi_1(M)]=\infty$ for $k\gg 1$, then a finite cover $M'$ of $M$ is a warped product (see Theorem \ref{thm:warped}), but no bound on the degree of the cover $M'\rightarrow M$ is obtained. The existence of a bound is proven in the following lemma.
	\begin{lemnr} There exists $d(n)\geq 1$ such that the cover $M'\rightarrow M$ can be chosen of degree at most $d(n)$.\label{lem:degbound}\end{lemnr}
	\begin{proof} We keep the notation of the proof of Proposition \ref{prop:prod}. By Lemma \ref{lem:fwclaims}.(\ref{lem:lattice}), we have that $\Gamma_0\subseteq I^0$ is a cocompact lattice and by Proposition \ref{prop:ss}, we know that $I^0$ is semisimple with finite center. As in the proof of Proposition \ref{prop:prod}, consider the the short exact sequence
		\begin{equation} 1\rightarrow I^0\rightarrow \langle I^0, \Gamma\rangle \rightarrow \Gamma\slash\Gamma_0\rightarrow 1,\label{eq:ses}\end{equation}
	which gives rise to morphism $\sigma:\Gamma\slash\Gamma_0\rightarrow\Out(I^0)$. Let $\Gamma'$ be the preimage in $\langle I^0, \Gamma\rangle$ of $\ker\sigma$. Now consider the short exact sequence
		\begin{equation}
			1\rightarrow I^0\slash Z(I^0)\rightarrow \langle \Gamma', I^0\rangle\slash Z(I^0)\rightarrow \Gamma'\slash\Gamma_0\rightarrow 1\label{eq:ses2}.
		\end{equation}
	This extension is determined by a cohomology class in $H^2(\Gamma'\slash\Gamma_0, Z(I^0\slash Z(I^0)))$ and a morphism $\Gamma'\slash\Gamma_0\rightarrow \Out(I^0\slash Z(I^0))$. Since $Z(I^0\slash Z(I^0))=1$ and $\Gamma'\slash\Gamma_0= \ker\sigma$, we see this extension is trivial, so that
		$$\langle \Gamma', I^0\rangle\slash Z(I^0)\cong (I^0\slash Z(I^0)) \times (\Gamma'\slash\Gamma_0).$$
	In particular $\Gamma'\slash\Gamma_0$ centralizes the image of $\Gamma_0$ in $I^0\slash Z(I^0)$. Since $Z(I^0)$ is finite and $\Gamma_0$ is torsion-free, it follows that $\Gamma_0$ projects isomorphically into $I^0\slash Z(I^0)$. Hence $\Gamma'\slash\Gamma_0$ centralizes $\Gamma_0$. 
	
	Further $Z(\Gamma_0)=\Gamma_0\cap Z(I^0)=1$, so we have $\Gamma'\cong \Gamma_0\times (\Gamma'\slash\Gamma_0)$. Now it is clear that 
		$$[\Gamma:\Gamma']\leq |\Out(I^0)|\leq |\Out(I^0\slash Z(I^0))|.$$
	To have this bound only depend on $n$ (not on $I^0$), let $d$ be the maximal order of $\Out(G)$ where $G$ is a semisimple Lie group with trivial center such that the associated symmetric space has dimension $\leq n$.\end{proof}	

We will now show that $[I(\widetilde{M},g_k):\pi_1(M)]=\infty$ for $k\gg 1$, which will complete the proof. By Theorem \ref{thm:coverbddbound}, it follows that there exists a nondiscrete Lie group $G$ (possibly with infinitely many components) such that $\Gamma\subseteq G$, and $G$ acts properly by $C^2$ diffeomorphisms on $\widetilde{M}$. Set $\Gamma_0:=\Gamma\cap G^0$. 

By Lemma \ref{lem:fwclaims} and Proposition \ref{prop:ss}, we have that $G^0$ is semisimple with finite center and no compact factors, and $\Gamma_0\subseteq G^0$ is a cocompact lattice. By Lemma \ref{lem:degbound} we can find $\Gamma'\subseteq \Gamma$ of index at most $d$ such that
	$$\Gamma'\cong \Gamma_0\times (\Gamma'\slash\Gamma_0).$$
By Theorem \ref{thm:eghconv} of Fukaya-Yamaguchi, we find that there is a closed subgroup $H\subseteq G$ such that 
	$$(\widetilde{M}, I_k, p_k)\to(\widetilde{M}, H, p_k)$$
in the equivariant Gromov-Hausdorff sense. Note that since $H\subseteq G$ is closed, $H$ is a Lie group.

\textbf{Step 1 (Structure of the limit)}. We claim that $H$ is a nontrivial product of factors of $G^0$ and $\Lambda_0:=\Gamma_0\cap H^0$ is a cocompact lattice in $H^0$.

The proof of Theorem \ref{thm:eghconv} (see \cite[Prop 3.6]{fyalmostneg}) shows that if $f_k\in I_k$ such that $f_k\rightarrow f$ (or along a subsequence) then $f\in H$. It follows that $\Gamma\subseteq H$ and from the proof of Theorem \ref{thm:coverbddbound}, it follows that there are infinitely many distinct cosets of $G\slash\Gamma$ that are limits of subsequences of $(f_k)_k$ where $f_k\in I_k$. Hence $[H:\Gamma]=\infty$. 

Since $H$ is a Lie group and $\Gamma\subseteq H$ is cocompact and of infinite index, it follows that $H^0\neq 1$. Further $\Gamma$ normalizes $H^0$. Hence $H^0$ is a closed connected subgroup of the semisimple Lie group $G^0$ that is normalized by $\Gamma_0$. By the Borel Density theorem (see \cite[5.17, 5.18]{raghlie}), it follows that $H^0$ is normal in $G^0$, hence $H^0$ is a product of factors of $G^0$. Since $\Gamma$ is cocompact in $H$, it follows that $\Lambda_0$ is cocompact in $H^0$. This completes Step 1.

\textbf{Step 2 (Detection of $H^0$ in the sequence)}. We show that the conditions for Theorem \ref{thm:nordetect} are satisfied. Namely, we claim 
	\begin{enumerate}
		\item $H\slash H^0$ is discrete and finitely presented,
		\item $\widetilde{M}\slash G$ is compact,
		\item There exists $R>0$ such that $H^0$ is generated by $H^0(R)$,
		\item $I_k$ acts on $\widetilde{M}$ properly discontinuously.
	\end{enumerate}
\emph{Proofs}.
	\begin{enumerate}
		\item It is clear that $H\slash H^0$ is discrete, and using that $\Gamma\cong\Gamma_0\times (\Gamma\slash\Gamma_0)$, it follows that
		$$H\slash H^0\cong (\Gamma_0\slash\Lambda_0)\times (\Gamma\slash\Gamma_0).$$
		$\Gamma$ is finitely presented since it is the fundamental group of a closed manifold. Note that $\Gamma\slash\Gamma_0$ is finitely presented since it is a direct factor of $\Gamma$ (see \cite[Lemma 1.3]{wallfin}). Further $\Gamma_0\slash\Lambda_0$ is finitely presented since it is a cocompact lattice in the semisimple Lie group $G^0\slash H^0$. Hence $H\slash H^0$ is finitely presented.
		\item Since $M$ is compact and the map
		$$M=\widetilde{M}\slash\Gamma\rightarrow \widetilde{M}\slash G$$
		is continuous and surjective, it follows that $\widetilde{M}\slash G$ is compact. 
		\item For any $R>0$, we know that $H^0(R)$ is an open neighborhood of the identity. Since $H^0$ is a connected Lie group, it is generated by any open neighborhood of the identity. 
		\item The action of $I_k$ on $\widetilde{M}$ is properly discontinuous since $I_k$ contains $\Gamma$ with finite index, and $\Gamma$ acts properly discontinuously on $\widetilde{M}$.
	\end{enumerate}
By Properties (1)-(4) the hypotheses of Theorem \ref{thm:nordetect} hold, so there exist normal subgroups $I_k'\subseteq I_k$ such that $I_k\slash I_k'\cong H\slash H^0$ for $k\gg 1$. 

\textbf{Step 3 (Constructing maps to $H^0$)}. Since $H\slash H^0\cong \Gamma\slash\Lambda_0$, it follows that $I_k'$ contains $\Lambda_0$ with finite index. Choose $\Lambda_k\subseteq\Lambda_0$ of finite index such that $\Lambda_k$ is normal in $I_k'$. Therefore we get a map
	$$I_k'\slash\Lambda_k\rightarrow \textrm{Out}(\Lambda_k).$$
The generalized Nielsen realization problem asks if whenever $N$ is a closed aspherical manifold and $F\subseteq \Out(\pi_1(N))$ is a finite subgroup, $F$ can be realized as a group of isometries $N$. Originally Nielsen posed this problem for hyperbolic surfaces and asked if $F$ can be realized by isometries. An affirmative solution was given by Kerckhoff \cite{nielsenreal}. For locally symmetric spaces of noncompact type without surface factors Mostow rigidity implies that $F$ lifts to a group of isometries \cite{mostowrigid}. A combination of these results solves the generalized Nielsen realization problem for arbitrary closed locally symmetric spaces of noncompact type.

It follows that $I_k'\slash\Lambda_k$ acts by isometries on the locally symmetric space $\Lambda_k \backslash H^0\slash K$, where $K\subseteq H^0$ is a maximal compact subgroup. It follows that there is a map 
	$$\varphi_k: I_k'\rightarrow \isom(H^0\slash K).$$
This map extends the natural inclusion $\Lambda_0\hookrightarrow H^0$. The number of components of $\isom(H^0\slash K)$ is finite, so that there exists $L>1$ (independent of $k$) and a characteristic subgroup $I_k''\subseteq I_k'$ of index at most $L$ such that
	$$\varphi_k: I_k''\rightarrow H^0.$$
\textbf{Step 4 (Estimating $\ker\varphi_k$)}. We show that $\ker\varphi_k\neq 1$ for $k\gg 1$. By the Kazhdan-Margulis theorem on minimal orbifolds for symmetric spaces \cite{minorbhomog}, the volume of orbifolds modeled on $H^0\slash K$ is bounded away from 0. Since
	$$[\varphi_k(I_k''):\Lambda_0]=\frac{\vol(\Lambda_0\backslash H^0\slash K)}{\vol(\varphi_k(I_k'')\backslash H^0\slash K)},$$
it follows that there is $L'$ (independent of $k$) such that $[\varphi_k(I_k''):\Lambda_0]\leq L'$. On the other hand we know
	$$[I_k'':\Lambda_0]\geq \frac{1}{L}[I_k':\Lambda_0]\rightarrow\infty.$$
Hence we have
	$$|\ker\varphi_k|\rightarrow\infty.$$
\textbf{Step 5 (End of the proof)}. Let $\overline{F}_k$ be the center of $\varphi_k(I_k'')$, and set $F_k:=\varphi_k^{-1}(\overline{F}_k).$ Since $\varphi_k(I_k'')$ is a lattice in the connected semisimple Lie group $H^0$ (with finitely many components and finite center), it follows from the Borel Density theorem that any finite normal subgroup of $\varphi_k(I_k'')$ is central (see \cite[Cor 4.45]{wmarithm}). Hence $F_k$ is characteristic in $I_k''$. Since $I_k''$ is normal in $I_k$, it follows that $F_k$ is normal in $I_k$. 

Now let $\Gamma'\subseteq \Gamma$ be the centralizer of $F_k$. Then $\Gamma'\subseteq\Gamma$ is normal and finite index. Hence the action of $F_k$ on $\widetilde{M}$ descends to an action on the closed aspherical manifold $M':=\widetilde{M}\slash\Gamma'$.
 
 By assumption, $\Gamma$ has no nontrivial normal abelian subgroups. On the other hand $Z(\Gamma')$ is a characteristic subgroup of $\Gamma'$, hence $Z(\Gamma')$ is a normal abelian subgroup of $\Gamma$. Therefore $Z(\Gamma')=1$.
 
 By a theorem of Borel \cite{borelnoact}, if $F$ is a finite group acting effectively on a closed aspherical manifold with centerless fundamental group $\pi$, then the induced map $F\rightarrow \Out(\pi)$ is injective. But $F_k$ centralizes $\Gamma'$, so the map $F_k\rightarrow \Out(\Gamma')$ is trivial. We conclude that $F_k$ is trivial. But since $F_k$ contains $\ker\varphi_k$, we know that $|F_k|\rightarrow\infty$. This is a contradiction.\end{proof}

\section{Minimal Orbifolds}
\label{sec:minorb}

\begin{proof}[Proof of Corollary \ref{cor:minorbgen}] The proof is by induction on $n$. First note that $\injrad(M)\geq\varepsilon$ combined with Berger's isembolic inequality \cite{isemineq} yields
	$$\vol(M)\geq v$$
for some $v>0$ only depending on $n$ and $\varepsilon$. Let $C$ be as in Theorem \ref{thm:coverreg}. First suppose we have $[I(\widetilde{M}):\pi_1(M)]\leq C$. Then
	\begin{align*}
	\vol(\widetilde{M}\slash I(\widetilde{M}))&=\frac{\vol(M)}{[I(\widetilde{M}):\pi_1(M)]}\\
											  &\geq \frac{\vol(M)}{C}\\
											  &\geq \frac{v}{C},
	\end{align*}
which yields the desired bound. 

Now suppose we have $[I(\widetilde{M}):\pi_1(M)]>C.$ By Theorem \ref{thm:coverreg} and choice of $C$, we have there exists a contractible Riemannian manifold $X$ and a nontrivial symmetric space $Y$ of noncompact type such that $\widetilde{M}$ is isometric to the Riemannian warped product $X\times_f Y$ where $f:X\rightarrow \bbR_{>0}$ is a smooth map. Therefore at the point $(x,y)\in X\times Y$ the metric $g$ satisfies
	$$g|_{(x,y)}=g_X|_x \oplus f(y) g_Y|_y$$
where $g_X$ is a metric on $X$ and $g_Y$ is a locally symmetric metric on $Y$. 

\subsubsection*{Base case ($n=2$)} Since the only nontrivial locally symmetric space of dimension $\leq 2$ is the hyperbolic plane, we must have $Y\cong \bbH^2$, and $X$ is trivial. In particular it follows that $M$ is a hyperbolic surface of constant curvature $\geq -\Lambda$. Siegel \cite{siegeldisc} showed that the minimal Euler characteristic of a 2-orbifold is $\frac{1}{42}$, the Gauss-Bonnet theorem implies that for any orbifold $M'$ modeled on $\widetilde{M}$, we have
	$$\vol(M)\geq \frac{-2\pi\chi(M)}{\Lambda}\geq \frac{\pi}{21\Lambda}.$$
	
\subsubsection*{Inductive step} Suppose now the result is true in dimension $<n$. Normalize $g_Y$ such that $K(g_Y)\geq -1$.

\textbf{Step 1 (Control on geometry in the $X$-direction)}. Let $d\geq 1$ be as in Theorem \ref{thm:coverreg}. Choose a cover $M'$ of degree $\leq d$ that is isometric to a warped product $B\times_f N$ where $N$ is a locally symmetric space modeled on $Y$. We can write $N=Y\slash\Lambda_0$ and $B=X\slash(\Lambda\slash\Lambda_0)$ for some $\Lambda\subseteq I(\widetilde{M})$. Since the factor $X$ is totally geodesic in the warped product $X\times_f Y$ and the cover $M'\rightarrow M$ is degree $\leq d$, we see that
	$$|\Ric_{g_X}|\leq \Lambda, \hspace{0.5 cm} \injrad(B)\geq\varepsilon, \hspace{0.5 cm} \diam(B)\leq dD.$$
Further $B$ is aspherical and $\Lambda\slash\Lambda_0$ does not contain nontrivial normal abelian subgroups, and $\dim B<n$. By the inductive hypothesis there exists $\mu'>0$ (only depending on $n, \Lambda, \varepsilon, D$) such that
	$$\vol(X\slash(I\slash I^0))\geq \mu',$$
where $I:=\isom(\widetilde{M})$.

\textbf{Step 2 (Control in $Y$-direction)} We claim there exists $\alpha>0$ only depending on $\Lambda, \varepsilon, D$ and $n$, such that $f\geq \alpha$ everywhere.

To see this, choose $x_0\in X$ such that $f(x_0)$ is a minimum. Then $\{x_0\}\times Y$ is totally geodesic in $X\times_f Y$. This is immediate from an explicit description of the geodesics of a warped product (see \cite[Theorem 6.3]{zeghibwarped}). Namely, suppose $\gamma(t)=(x(t),y(t))\in X\times_f Y$ is a geodesic. Then $y(t)$ is an unparametrized geodesic on $Y$, and $x(t)$ satisfies the equation
	$$x''=-\nabla\left(\frac{c}{f}\right)(x)$$
for some $c>0$. Since $x(0)=x_0$ is a minimum of $f$, we see that $-\nabla\left(\frac{c}{f}\right)(x_0)=0$, so that $x(t)=0$ is a solution.

Since $\{x_0\}\times Y$ is totally geodesic in $X\times_f Y$, we know that the metric $f(x_0)g_Y$ also satisfies the Ricci curvature bound $|\Ric_{f(x_0)g_Y}|\leq\Lambda$. The normalization $K(g_Y)\geq -1$ and the bound $\Ric(f(x_0)g_Y)\geq -\Lambda$ give a lower bound $f(x_0)\geq \alpha$ where $\alpha$ depends only on $\Lambda$ and $n$. This completes the proof of Step 1.

Now suppose $\Delta$ acts on $\widetilde{M}$ properly discontinuously and $\vol(\widetilde{M}\slash\Delta)<\infty$. Set $\Delta_0:=\Delta\cap I^0$. Then $\Delta_0$ is a lattice in $I^0$. In the cocompact case this follows from \cite{FW}. For the noncocompact case, see the proof of \cite[Corollary 1.5]{localriemsymm}. To summarize, disintegrate the finite volume measure on $\widetilde{M}\slash\Delta$ along the fibers of the map 
	$$p:\widetilde{M}\slash\Delta\rightarrow X\slash(\Delta\slash\Delta_0).$$
If $x$ is not a singular point of the orbifold $X\slash(\Delta\slash\Delta_0)$, then $p^{-1}(x)=Y\slash\Delta_0$ equipped with the volume form $\nu_x:=f(x)^m\vol_Y$ where $m:=\dim Y$. Since the total space $\widetilde{M}\slash\Delta$ has finite volume and $f$ is bounded away from 0, it follows that a.e. fiber has finite volume.

By the Kazhdan-Margulis theorem on minimal orbifolds of symmetric spaces, we know there exists $\eta>0$ such that
	$$\vol_Y(Y\slash\Delta_0)\geq\eta.$$
\textbf{Step 3 (Bounding volume of orbifolds)}. We conclude that
	\begin{align*}\vol(\widetilde{M}\slash\Delta)&=\int_{X\slash(\Delta\slash\Delta_0)} \nu_x(p^{-1}(x)) \dvol_X(x)\\
												 &=\int_{X\slash(\Delta\slash\Delta_0)} |f(x)|^m \vol_Y(Y\slash\Delta_0)\dvol_X(x)\\
												 &\geq \alpha^m \vol_Y(Y\slash\Delta_0)\vol_X(X\slash(\Delta\slash\Delta_0))\\
												 &\geq \alpha^n\eta\mu'.
	\end{align*}
This proves the theorem.\end{proof}

\section{Detecting normal subgroups and equivariant Gromov-Hausdorff convergence}
\label{sec:nordetect}

We prove Theorem \ref{thm:nordetect}. As mentioned in Section \ref{sec:prelim}, the proof is essentially Fukaya-Yamaguchi's proof of Theorem \ref{thm:fynordetect} (see \cite[Appendix 1]{fyalmostneg}), but some arguments involving covering spaces are replaced by analogous arguments involving orbifold covers. Since the original proof is quite long, we will merely summarize most of the proof and only supply the details where the arguments need to be changed (see Claims \ref{cl:orbicover} and \ref{cl:orbihomeo} below). For the details of the original proof, see \cite{fyalmostneg}.
\begin{proof}[Proof of Theorem \ref{thm:nordetect}] Let $\varepsilon_k\rightarrow 0$ such that there are $\varepsilon_k$-equivariant Gromov-Hausdorff approximations
	$$(f_k, \varphi_k, \psi_k): (X_k, \Gamma_k, p_k)\rightarrow (Y, G, q).$$
Choose $R$ such that $R>R_0$ and $R>D$ and $R<\frac{1}{10\varepsilon_k}$. Then we define
	$$\Gamma_k'(R):=\Gamma_k(R) \cap \varphi_k^{-1}(G').$$
Let $\Gamma_k''$ be the subgroup generated by $\Gamma_k'(R)$. We will construct $\Gamma_k'\lhd \Gamma_k$ such that $\Gamma_k'\cap \Gamma_k(R)=\Gamma_k'(R)$ ,and $\Gamma_k\slash\Gamma_k' \cong G\slash G'$ for $k$ sufficiently large.

First define a relation on $\Gamma_k(R)$ as
	$$\gamma\sim \delta \Longleftrightarrow \gamma^{-1}\delta\in \Gamma_k'(3R).$$
As in \cite{fyalmostneg}, we show $\sim$ is a bi-invariant equivalence relation if $k$ is sufficiently large. Set
	$$\Lambda_k(R):=\Gamma_k(R)\slash\sim.$$
Since $\Gamma_k(R)$ is naturally a pseudogroup, there is a natural pseudogroup structure on $\Lambda_k(R)$. We apply the same construction to $G(R)$ and $G'(R)$ (rather than $\Gamma_k(R)$ and $\Gamma_k'(R)$), and we obtain a pseudogroup $H(R):=G(R)\slash\sim$. The maps $\varphi_k$ descend for $k\gg 1$ to pseudogroup morphisms
	$$\overline{\varphi}_k:\Lambda_k(R)\rightarrow H(R),$$
and for $k\gg 1$, we see that $\overline{\varphi}$ is an isomorphism. We let $\widehat{\Lambda_k(R)}$ and $\widehat{H(R)}$ be the groupifications of $\Lambda_k(R)$ and $H(R)$. There are natural inclusions
	$$i_k:\Lambda_k(R)\hookrightarrow \widehat{\Lambda_k(R)}$$
and
	$$i:H(R)\hookrightarrow \widehat{H(R)}.$$
Now $\Gamma_k'(3R)$ (resp. $G'(3R)$) acts (as a pseudogroup) on $B_{X_k}(p_k;R)$ (resp. $B_Y(q;R)$). We let $V_k(R)$ (resp. $W(R)$) be the quotients. 
\begin{claimnr} The natural projection $\pi: B_{X_k}(p_k;R)\rightarrow V_k(R)$ is an orbifold covering.\label{cl:orbicover}\end{claimnr}
	\begin{proof} Let $x\in B_{X_k}(p_k;R)$. Choose an open neighborhood $U\ni x$ such that $U\subseteq B_{X_k}(p_k;R)$ and $U$ contains no $\Gamma$-translates of $x$ except $x$ itself. On $U$, we have the local description
		$$U\rightarrow U\slash \textrm{Stab}_{\Gamma_k'(3R)}(x)$$
	for $\pi$. Therefore it suffices to prove $\textrm{Stab}_{\Gamma_k'(3R)}(x)$ is a group (note that this is not obvious, since $\Gamma_k'(3R)$ is only a pseudogroup). Let $\gamma,\delta\in \textrm{Stab}_{\Gamma_k'(3R)}(x)$. We want to show that $\gamma\delta^{-1}\in \Gamma_k'(3R)$. First we show that $\gamma\delta^{-1}\in \Gamma_k(3R)$. Note that
		\begin{align*}
		d_k(\gamma\delta^{-1} p_k, p_k)&\leq d_k(\gamma\delta^{-1}, \gamma\delta^{-1}x) + d_k(\gamma\delta^{-1}x,x)+d_k(x,p_k)\\
									&	\leq R+0+R\\
									&\leq 2R,
		\end{align*}
	so $\gamma\delta^{-1}\in \Gamma_k(3R)$. It remains to show that $\gamma\delta^{-1}\in \varphi_k^{-1}(G')$. To see this, note that since $\overline{\varphi}_k$ is a pseudogroup morphism, we have
		$$\overline{\varphi}_k([\gamma\delta^{-1}])=[e]\in H(R),$$
	so $\varphi_k(\gamma\delta^{-1})\in G'$, as desired. It follows that $V_k(R)$ is an orbifold and $B_{X_k}(p_k;R)\rightarrow V_k(R)$ is an orbifold cover.\end{proof}

The quotient pseudogroup $\Lambda_k(3R)$ (resp. $H(3R)$) acts on $V_k(R)$ (resp. $W(R)$). Hence $\Lambda_k(3R)$ (resp. $H(3R)$) acts (as a pseudogroup) on $\widehat{\Lambda_k(3R)}\times V_k(R)$ (resp. $\widehat{H(3R)}\times W(R)$) as
	$$\gamma\cdot (\delta, x)=(\delta\gamma^{-1}, \gamma x)$$
for $\gamma\in \Lambda_k(3R), \delta\in\widehat{\Lambda_k(3R)}$ and $x\in V_k(R)$. There is a similar formula for $H(3R)$ acting on $\widehat{H(3R)}\times W(R)$. 
Let $V_k'$ (resp. $W'$) be the quotient, and let $V_k$ (resp. $W$) be the connected component of the image of $[(e,[p_k])]\in V_k'$ (resp. $[(e,[q])]$).

Since the action of $\Lambda_k(3R)$ on $\widehat{\Lambda_k(3R)}\times V_k(R)$ is free, we have a natural orbifold structure on $V_k$. 

Further $\widehat{\Lambda_k(3R)}$ (resp. $\widehat{H(3R)}$) acts on $\widehat{\Lambda_k(3R)}\times V_k(R)$ (resp. $\widehat{H(3R)}\times W(R)$ by left-translations on the first factor, and this action commutes with the action of $\Lambda_k(3R)$ (resp. $H(3R)$) described above. Therefore $\widehat{\Lambda_k(3R)}$ (resp. $\widehat{H(3R)}$) acts on $V_k'$ (resp. $W'$). Let $\Lambda_k$ (resp. $H$) be the subgroups preserving $V_k$ (resp. $W$). 

As in Fukaya-Yamaguchi's proof, we see:
	\begin{enumerate}
		\item $\Lambda_k\cong H$,
		\item $\Lambda_k$ (resp. $H$) acts properly discontinuously on $V_k$ (resp. $W$),
		\item $W\slash H\cong Y\slash G$.
	\end{enumerate}
\begin{claimnr}[replaces Lemma A1.14 in \cite{fyalmostneg}] $V_k\slash\Lambda_k\cong X_k\slash\Gamma_k$ as orbifolds.\label{cl:orbihomeo}\end{claimnr}
\begin{rmknr} Lemma A1.14 of \cite{fyalmostneg} states that $V_k\slash \Lambda_k$ and $X_k\slash\Gamma_k$ are diffeomorphic. Here it is critical that $\Gamma_k$ act freely, for otherwise $X_k\slash\Gamma_k$ need not be a manifold. Below we will give an explicit description of orbifold charts which will imply the result.\end{rmknr}
\begin{proof}Since $(X_k, \Gamma_k, p_k)\to(Y, G, q)$, we have Gromov-Hausdorff convergence of the orbit spaces (see Theorem \ref{thm:orbitgh}):
	$$X_k\slash\Gamma_k\to Y\slash G.$$
Hence for $k\gg 1$ we have $\diam(X_k\slash\Gamma_k)<R$. It follows that we can choose a fundamental domain $U_k\subseteq B_{X_k}(p_k;R)$ for the $\Gamma_k$-action on $X_k$. 

Then $U_k$ projects to the fundamental domain $F_k:=U_k\slash\Gamma_k'(3R)$ for the $\Lambda_k$-action on $V_k$. If $x\in F_k\subseteq V_k$ and $g\in\Lambda_k$ such that $g x = x$, then there is $\gamma\in\Lambda_k(3R)$ such that
	$$(g,x)=(\gamma^{-1},\gamma x).$$
Hence $g=\gamma^{-1}$, so $g\in \textrm{Stab}_{\Lambda_k(3R)}(x)$ (here we view $x\in U_k\slash\Gamma_k'(3R)\subseteq V_k(R)$). Choose a lift $\widetilde{x}\in U_k$ of $x$. Then a chart for $V_k\slash\Lambda_k$ containing $x$ is given by $(U_k, \textrm{Stab}_{\Gamma_k(3R)}(\widetilde{x}), \pi_1)$ where $\pi_1$ is the natural projection $\pi: U_k\rightarrow V_k\slash U_k\slash \textrm{Stab}_{\Gamma_k(3R)}(\widetilde{x})$. Again an important point is that $\textrm{Stab}_{\Gamma_k(3R)}(\widetilde{x})$ is a group, even though $\Gamma_k(3R)$ is only a pseudogroup. This is proven in the same way as in Claim \ref{cl:orbicover}.

A chart for $X_k\slash\Gamma_k$ containing an orbit $x=\Gamma_k \widetilde{x}$ is given by the triple $(U_k, \textrm{Stab}_{\Gamma_k}(\widetilde{x}), \pi_2)$ where $\pi_2$ is the natural projection $U_k\rightarrow X_k\slash\Gamma_k$. Since
	$$\textrm{Stab}_{\Gamma_k}(\widetilde{x})=\textrm{Stab}_{\Gamma_k(3R)}(\widetilde{x}),$$
it follows that $V_k\slash\Lambda_k\cong X_k\slash\Gamma_k$ as orbifolds.\end{proof}
In particular, $V_k\slash\Lambda_k$ is a good orbifold (because $X_k$ is a manifold) and
	$$V_k\rightarrow V_k\slash \Lambda_k$$
is an orbifold covering. Therefore $V_k$ is a good orbifold as well. Let $\Gamma_k'$ be the orbifold fundamental group of $V_k$. With this setup, the remainder of Fukaya-Yamaguchi's proof works verbatim.
\end{proof}

\bibliographystyle{alpha}
\bibliography{localsymmbib}

\end{document}